\documentclass[12pt,]{amsart}
\usepackage[latin1]{inputenc}
\usepackage{indentfirst}
\usepackage[dvips]{graphicx}
\usepackage{amsfonts}
\usepackage{amssymb}
\usepackage{theoremref}
\usepackage{amsthm}
\usepackage{bbold}
\usepackage{dsfont}
\usepackage{amsmath}
\usepackage{pb-diagram}
\usepackage{geometry}
\usepackage{pstricks,pst-node}
\usepackage{lmodern}
\usepackage{etoolbox}
\apptocmd{\sloppy}{\hbadness 10000\relax}{}{}
\newcommand {\nl}{\newline}

\begin{document}
\newtheorem{defi}{Definition}[section]
\newtheorem{prop}[defi]{Proposition}
\newtheorem{teo}[defi]{Theorem}
\newtheorem{lemma}[defi]{Lemma}
\newtheorem{corol}[defi]{Corollary}
\newtheorem{obs}[defi]{Remark}
\newtheorem{ex}[defi]{Example}
\newcommand {\dem}{\textbf{Proof: }}
\newcommand {\afi}{\textbf{Claim: }}
\newcommand {\prafi}{\textbf{Proof of claim: } }
\newcommand {\D}{\displaystyle}
\newcommand {\eps}{\varepsilon}\title{C$^*$-algebras associated with endomorphisms of groups}

\author{F. Vieira}\thanks{Supported by CAPES - Coordena\c{c}\~{a}o de Aperfei\c{c}oamento de Pessoal de N\'{i}vel Superior}
\date{}

\begin{abstract}
In this work we construct a C$^*$-algebra from an injective
endomorphisms of some group $G$, allowing the endomorphism to have
infinite cokernel. We generalize results obtained by
I. Hirshberg in \cite{Hirsh} and by J. Cuntz and A. Vershik in \cite{CunVer}. In good cases we show that
the C$^*$-algebra that we study is classifiable by Kirchberg's
classification theorem, with K-groups equal to $K_{*}(C^*(G))$.
\end{abstract}
\maketitle

\vspace{1cm}
\section{Introduction}

In \cite{Hirsh} Hirshberg defined a C$^*$-algebra associated with endomorphisms of groups with finite cokernel. The obvious sequence of that paper is to construct the same C$^*$-algebra for endomorphisms with infinite cokernel. So in this paper we define and study a universal C$^*$-algebra constructed from an injective endomorphism $\varphi$ with infinite cokernel of a discrete countable group. Thus the biggest difference of this paper with Hirshberg's (\cite{Hirsh}) is that we allow
\begin{equation*}
\left|\dfrac{G}{\varphi(G)} \right|=\infty.
\end{equation*}

In order to generalize the constructions, we also associate the C$^*$-algebra with a set $B$ of subgroups of $G$ and call it $\mathds{U}[\varphi,B]$. Their r\^{o}le is to implement naturally the multiplication rule inside $\mathds{U}[\varphi,B]$, because here we do not have finitely many projections summing up to one. The relations defining $\mathds{U}[\varphi,B]$ are dictated by the natural representations of $\varphi$, $B$ and $G$ on the Hilbert space $l^2(G)$ of all square summable complex functions on $G$. The unitaries representing the group elements, the projections associated with subgroups of $G$ and the isometry representing $\varphi$ generate a concrete C*-subalgebra of $\mathcal{L}(l^2(G))$ which, in good cases, is isomorphic to the C*-algebra $\mathds{U}[\varphi,B]$ and can thus be described by generators and relations.

Beside Hirshberg's paper, similar constructions have been studied before by various authors \cite{Buss1}, \cite{Cuntz1}, \cite{Culi1}, \cite{CunVer}, \cite{Li1} and \cite{Li2}. In particular, somewhat similar C*-algebras have been associated with endomorphisms of abelian groups and also with semigroups. Also the ring C*-algebras studied in \cite{Cuntz1}, \cite{Culi1} arise in a similar way.

It is important to note that even though $\varphi(G)$ is only a subgroup (not necessarily normal, because $G$ can be non-abelian) we are able to count how many elements there are in the above quotient. Also choose some family $B$ of subgroups of $G$, containing $G$.

The group elements give rise to unitary operators $\{U_g\}_{g\in G}$ acting on $l^2(G)$ by left multiplication, and the endomorphism induces an isometry $S$ acting on $l^2(G)$ through $\varphi$: denoting by $\{\xi_h:\,h\in G\}$ the canonical orthonormal basis of $l^2(G)$, $S$ is defined by $S(\xi_h):=\xi_{\varphi(h)}$. With every element $H$ of $B$, consider the projection $E_{[{H}]}$ with
\begin{equation}\label{equa9}
E_{[H]}(\xi_g)=\left\{
                 \begin{array}{ll}
                   \xi_g, & \hbox{if }g\in H; \\
                   0, & \hbox{otherwise.}
                 \end{array}
               \right.
\end{equation}

The C$^*$-subalgebra of $\mathcal{L}(l^2(G))$ generated by the operators above is denoted $C_r^*[\varphi,B]$ (Definition \ref{defi2}).

The operators thus defined satisfy some natural relations, and we use these relations to define the universal C$^*$-algebra $\mathds{U}[\varphi,B]$, associated with $\varphi$. Particularly, the sum condition which appears in \cite{CunVer} does not hold in our situation (we would have a inequality of the form $<\infty$). However the projections associated with the subgroups of $G$ in $\mathds{U}[\varphi,B]$ take their place and are crucial to prove the main results.

One of the important ones obtained is that, if $G$ is amenable, if the intersection of the subgroups of $B$ contains some image of $G$ through $\varphi$ and if $\varphi$ is pure\footnote{The endomorphism $\varphi$ is pure when $\D\bigcap_{n\in\mathds{N}}\varphi^n(G)=\{e\}$.}, then $\mathds{U}[\varphi,B]$ is a Kirchberg algebra\footnote{$\mathds{U}[\varphi,B]$ is separable, nuclear, simple and purely infinite.}. In particular, this implies that in this case $\mathds{U}[\varphi]$ and $C_r^*[\varphi]$ are isomorphic. This result also extends the ones obtained by Hirshberg in \cite{Hirsh} and by Cuntz and Vershik in \cite{CunVer}.

To prove the result above, it is crucial to use a semigroup crossed product description of $\mathds{U}[\varphi,B]$. Here we use the definition of a semigroup crossed product presented by Li in Appendix A of \cite{Li1} using covariant representations. The semigroup implementing the crossed product can be the semidirect product $S:=G\rtimes_\varphi\mathds{N}$ or the semigroup of natural numbers $\mathds{N}$. Such a description allows us to use the six term exact sequence presented by Khoshkam and Skandalis \cite{Khoska} to calculate the K-theory of our C*-algebra, as is done by Cuntz and Vershik in \cite{CunVer}.

The semigroup crossed product description above also implies the existence of a (full corner) group crossed product description of $\mathds{U}[\varphi,B]$ (by the minimal automorphic dilation introduced in \cite{Cuntz2}, \cite{CuntzTopMarkovII} and generalized later by M. Laca in \cite{Laca1}), using the group of integers $\mathds{Z}$. This allows one to use the classical Pimsner-Voiculescu exact sequence \cite{Pivo1} to calculate their K-groups. 

We will see that considering $B=\{G\}$ gives interesting examples, and we then denote the C$^*$-algebra only by $\mathds{U}[\varphi]$. In this case, the isomorphism above is not the only way to represent it as a crossed product: analogously to the work of G. Boava and R. Exel in \cite{BoEx} one can show that $\mathds{U}[\varphi]$ has a partial group crossed product description, which can also be related to an inverse semigroup crossed product by \cite{ExVi}. Apart from giving $\mathds{U}[\varphi]$ another description by an established structure, this result also provides another way to prove the simplicity of $\mathds{U}[\varphi]$ in some cases.

It can be noted that a particular semigroup is very important in our constructions: the semigroup $S=G\rtimes_\varphi\mathds{N}$. We prove that when the group $G$ is amenable and $\varphi$ is pure, the three semigroup C$^*$-algebras defined by Li in \cite{Li2} - namely C$^*(S)$, C$^*_s(S)$ and C$^*_r(S)$ - associated with the semigroup $S$ are isomorphic to $\mathds{U}[\varphi]$ and also nuclear, simple and purely infinite (Theorem \ref{teo4}), answering partially one open question in \cite{Li2}.

To finish, using the semigroup crossed product description of $\mathds{U}[\varphi]$ from Chapter 2, we study its K-theory. Using a natural split exact sequence and the six term exact sequence provided by Khoshkam and Skandalis \cite{Khoska} we easily conclude that the K-groups of $\mathds{U}[\varphi]$ are the same as the ones of $C^*(G)$. This implies that, imposing some extra conditions, $\mathds{U}[\varphi]$ is classifiable\footnote{$\mathds{U}[\varphi]$ is a Kirchberg algebra which satisfies UCT and $[1]_0=1$} by Kirchberg's classification theorem \cite{Kirchcla}.

A weaker version of the result concerning the K-theory of $\mathds{U}[\varphi]$ can be obtained independently using some recent results by Cuntz, Echterhoff and Li in \cite{CuEcLi1}.

I would like to thank J. Cuntz for the Ph.D. orientation and the helpful comments and corrections about this paper.

\section{Definitions and basic results}

In this paper $G$ will always be a discrete countable group with unit $e$ and $\varphi$ an injective endomorphism (monomorphism) of $G$ with infinite cokernel. When necessary, we require the amenability of $G$ or $\varphi$ to be pure. We want to construct a C$^*$-algebra associated with $\varphi$. To generalize Hirshberg's constructions even more, we also want to associate the C$^*$-algebra with some set $B$ of subgroups of $G$ which contains $G$. We consider it to have a natural behaviour of the multiplication rule inside the C$^*$-algebra.

Now we expand $B$ minimally to $C(B)$: so it is the smallest set of subsets of $G$ containing $B$ such that the \emph{regularity conditions} hold i.e, it is closed under finite unions, finite intersections, complements and under images of $\varphi$. Thus we have the following concrete C$^*$-algebra.
\begin{defi}\label{defi1}Consider $B$ a family of subgroups of $G$ (containing $G$) defined as above. We denote $C_r^*[\varphi,B]$ the reduced universal C$^*$-algebra generated by the three families of operators defined in (\ref{equa9}):
\begin{equation*}
\begin{split}
\hbox{a family of projections }&\{E_{[X]}:X\in C(B)\};\\
\hbox{unitaries }&\{U_g:g\in G\}\\
\hbox{and the isometry }&S.
\end{split}
\end{equation*}
\end{defi}

Studying the properties of the operators above, it is natural to define its universal version:
\begin{defi}\label{defi2}As above choose a set $B$ of subgroups of $G$ (containing $G$) and construct the family $C(B)$. Then $\mathds{U}[\varphi,B]$ is the universal C$^*$-algebra generated by
\begin{equation*}
\begin{split}
\hbox{a family of projections }&\{e_{[X]}:X\in C(B)\};\\
\hbox{unitaries }&\{u_g:g\in G\}\\
\hbox{and one isometry }&s
\end{split}
\end{equation*}

satisfying:
\begin{enumerate}
	\item[(i)] $u_gs^nu_hs^m=u_{g\varphi^n(h)}s^{n+m}$;
	\item[(ii)] $u_gs^ne_{[X]}{s^*}^nu_{g^{-1}}=e_{[g\varphi^n(X)]}$;
	\item[(iii)] $e_{[G]}=1$;
	\item[(iv)] $e_{[X]}e_{[Y]}=e_{[X\cap Y]}$ and
	\item[(v)] $e_{[X]}+e_{[Y]}=e_{[X\cup Y]}+e_{[X\cap Y]}$.
\end{enumerate}
\end{defi}

Since $u_gs^n{s^*}^nu_{g^{-1}}=e_{[g\varphi^n(G)]}$, the projections $u_gs^n{s^*}^nu_{g^{-1}}$ commute and considering $n\geq m$:
\begin{equation*}
\begin{split}
u_gs^n{s^*}^nu_{g^{-1}}u_hs^m{s^*}^mu_{h^{-1}}&=e_{[g\varphi^n(G)]}e_{[h\varphi^m(G)]}=e_{[g\varphi^n(G)\cap h\varphi^m(G)]}=\\
&=\nl \left\{
    \begin{array}{ll}
      e_{[g\varphi^n(G)]}, & \hbox{if }h\in g\varphi^m(G); \\
      0, & \hbox{otherwise,}
    \end{array}
  \right.\\
&=\nl \left\{
    \begin{array}{ll}
      u_gs^n{s^*}^nu_{g^{-1}}, & \hbox{if }h\in g\varphi^m(G); \\
      0, & \hbox{otherwise.}
    \end{array}
  \right.
\end{split}
\end{equation*}
\begin{obs}\em\label{obs1}It is important to mention that the construction above can be done when $B$ is just a set of subsets of $G$, or even consider some set $C$ in place of $C(B)$, containing any type of set and closed under the regularity conditions.
\end{obs}

First of all we will see that only the initial set $B$ is important to generate the C$^*$-algebras above, and the fact that in $C(B)$ some elements are not subgroups of $G$ is not a problem. Note that some elements of $C(B)$ are given by
\begin{equation*}
g\bigcap_{i=1}^{m}\varphi^{n_i}(H_i)
\end{equation*}

with $g\in G$, $n_i\in\mathds{N}$ and $H_i\in B$. In fact we can use these to describe the $*$-algebra $span(\{e_{[X]}:X\in C(B)\})$.
\begin{lemma}\label{lema1}Define
$$
B':=\left\{\bigcap_{i=1}^m\varphi^{n_i}(H_i): H_i\in B,\;n_i\in\mathds{N}\right\}.
$$

Then $span(\{e_{[X]}:X\in C(B)\})\cong span(\{e_{[gH']}:g\in G, H'\in B'\})=:D'$.
\end{lemma}
\dem $\supseteq:$ Obvious.

$\subseteq:$ Let us call $K':=\{X\subseteq G:e_{[X]}\in D'\}$. It is obvious that $B\subseteq K'$. Moreover $K'$ is closed under:
\begin{itemize}
\item $\D\bigcap_{i=1}^n:$ By definition, for $X_1,X_2\in K'$ it holds that
$$
e_{[X_1\cap X_2]}=e_{[X_1]}e_{[X_2]}.
$$

\item Complements: For $X\in K'$:
$$
e_{[X^c]}=1-e_{[X]}=e_{[G]}-e_{[X]}\in D'.
$$

\item $\D\bigcup^n_{i=1}:$ Note that $X\cup Y=[X^c\cap Y^c]^c\in K'$, $\forall\; X$, $Y\in K'$.\par\vspace{\baselineskip}

\item And if $X\in K'$, $g\in G$ and $n\in\mathds{N}$, the injectivity of $\varphi$ implies $g\varphi^n(X)\in K'$.\par\vspace{\baselineskip}
\end{itemize}

Therefore $K'$ satisfies the regularity conditions, and $C(B)\subseteq K'$, because $C(B)$ is the smallest set containing $B$ satisfying it. Then $span(\{e_{[X]}:X\in C(B)\})\subseteq D'$.
\begin{flushright}

  $\square$

  \end{flushright}

This result implies an important and simpler way to describe $\mathds{U}[\varphi,B]$:
\begin{prop}\label{prop1}The universal C$^*$-algebra $\mathds{U}[\varphi,B]$ is generated by
$$
\{e_{[H]}, u_g, s:H\in B, g\in G\}.
$$
\end{prop}
\dem Due to last lemma we only have to prove that
$$
span(\{e_{[gH']}:g\in G, H'\in B'\})\subseteq span\{e_{[H]}, u_g, s:H\in B, g\in G\}.
$$

But $e_{[gH']}=u_ge_{[H']}u_{g^{-1}}$ and for $H'=\D\bigcap_{i=1}^n\varphi^{n_i}(H_i)\in B'$ with $n_i\in\mathds{N}$ and $H_i\in B\cup\{G\}$ we have
$$
e_{[H']}=\prod^n_{i=1}e_{[\varphi^{n_i}(H_i)]}=\prod^n_{i=1}s^{n_i}e_{[H_i]}{s^*}^{n_i}.
$$

Therefore $e_{[gH']}\in C^*(\{e_{[H]}, u_g,s:H\in B, g\in G\})$.
\begin{flushright}

  $\square$

  \end{flushright}

\begin{obs}\em\label{obs2}Note that the lemma and the proposition above hold for any choice of $B$ (i.e, even if it doesn't consists of subgroups).
\end{obs}

Another interesting basic result:
\begin{prop}\label{prop2}Consider $\overline{B}$ containing only sets of the form $g_i\varphi^n(H_i)$, with $g_i\in G$ and $H_i$ subsets of $G$. Then
$$
\mathds{U}[\varphi,\overline{B}]\cong\mathds{U}[\varphi,B]
$$

where $B$ contains only the subsets $H_i$.
\end{prop}
\dem By Proposition \ref{prop1} (and the remark above),
\begin{equation*}
\begin{split}
\mathds{U}[\varphi,\overline{B}]&=C^*(\{e_{[g_iH_i]}, u_g, s:\;g_iH_i\in \overline{B}, g\in G\})\\
\mathds{U}[\varphi,B]&=C^*(\{e_{[H_i]}, u_g, s:\;H_i\in B, g\in G\})
\end{split}
\end{equation*}

But as
$$
e_{[g_i\varphi^n(H_i)]}=u_{g_i}s^ne_{[H_i]}{s^*}^nu_{g_i^{-1}}
$$

both C$^*$-algebras are isomorphic.
\begin{flushright}

  $\square$

  \end{flushright}

\begin{obs}\em\label{obs3}If we choose $B=\{G\}$ then $\mathds{U}[\varphi,B]$ is generated only by the unitary elements $\{u_g:\, g\in G\}$ and the isometry $s$, and it can be viewed as a natural generalization of the constructions in \cite{Hirsh} and \cite{CunVer}. This case will be studied with more details in Section \ref{seccan}.
\end{obs}

\section{Crossed product descriptions}

Define
$$
D[\varphi,B]:=C^*(\{u_gs^ne_{[H]}{s^*}^nu_{g^{-1}}:\, g\in G, n\in\mathds{N}, H\in B\})
$$

and note that it is a commutative C$^*$-subalgebra of $\mathds{U}[\varphi,B]$ because we have\nl $u_gs^ne_{[H]}{s^*}^nu_{g^{-1}}=e_{[g\varphi^n(H)]}$. We 
can define an action of the semigroup $S=G\rtimes_\varphi\mathds{N}$ on $D[\varphi,B]$ via
\begin{equation*}
\begin{split}
\alpha:S&\rightarrow \hbox{End}(D[\varphi,B])\\
   (g,n)&\mapsto u_gs^{n} (\cdot) {s^*}^nu_{g^{-1}}.
\end{split}
\end{equation*}

\begin{prop}\label{prop3}The C$^*$-algebra $\mathds{U}[\varphi,B]$ is isomorphic to $D[\varphi,B]\rtimes_{\alpha} S$.
\end{prop}
\dem In this proof we use the universality of both C$^*$-algebras to find the desired isomorphism. 
Remembering the definition from \cite{Li1}, $D[\varphi,B]\rtimes_{\alpha} S$ together with
\begin{equation*}
\begin{split}
\iota_D: D[\varphi,B]&\rightarrow D[\varphi,B]\rtimes_{\alpha} S\\
         x&\mapsto \iota_D(x)
\end{split}
\end{equation*}

and
\begin{equation*}
\begin{split}
\iota_S: S&\rightarrow \hbox{Isom}(D[\varphi,B]\rtimes_{\alpha} S)\\
  (g,n)&\mapsto \iota_S(g,n)
\end{split}
\end{equation*}

satisfying
$$
\iota_D(u_gs^nx{s^*}^nu_{g^{-1}})=\iota_S(g,n)\iota_D(x)\iota_S(g,n)^*
$$

is the semigroup crossed product of the dynamic system $(D[\varphi,B],S,\alpha)$.
But note that $\mathds{U}[\varphi,B]$ together with
\begin{equation*}
\begin{split}
\pi: D[\varphi,B]&\rightarrow \mathds{U}[\varphi,B]\\
         x&\mapsto x
\end{split}
\end{equation*}

and
\begin{equation*}
\begin{split}
\rho: S&\rightarrow \hbox{Isom}(\mathds{U}[\varphi,B])\\
  (g,n)&\mapsto u_gs^n
\end{split}
\end{equation*}

is a covariant representation of $(D[\varphi,B],S,\alpha)$, since:
$$
\rho(g,n)\pi(x)\rho(g,n)^*=u_gs^nx{s^*}^nu_{g^{-1}}=\pi(\alpha_{(g,n)}(x)).
$$

So we conclude that there exists a $*$-homomorphism
\begin{equation}\label{eq2}
\Phi: D[\varphi,B]\rtimes_{\alpha} S \rightarrow \mathds{U}[\varphi,B]
\end{equation}

such that $\Phi\circ\iota_D=\pi$ and $\Phi\circ\iota_S=\rho$.

In the other hand, it is well known that the crossed product $D[\varphi,B]\rtimes_{\alpha} S$ is generated as a C$^*$-algebra by elements of the form
$\iota_S(g,n)$ and $\iota_D(e_{[H]})$ with $H\in B$. Identifying $\iota_S(g,n)$ with
$u_gs^n$ and $\iota_D(e_{[H]})$ with $e_{[H]}$, it is easy to check that they satisfy conditions (i) - (v) of Definition \ref{defi2} which generate $\mathds{U}[\varphi,B]$.

Therefore we have a $*$-homomorphism
\begin{equation}\label{eq3}
\begin{split}
\Delta: \mathds{U}[\varphi,B] &\rightarrow D[\varphi,B]\rtimes_{\alpha} S\\
                      u_gs^n  &\mapsto \iota_S(g,n)\\
                     e_{[H]} &\mapsto \iota_D(e_{[H]}).
\end{split}
\end{equation}

We now show that (\ref{eq2}) and (\ref{eq3}) are inverses of each other:
$$
\Phi\circ\Delta(u_g)=\Phi(\iota_S(g,0))=\rho(g,0)=u_g
$$

$$
\Phi\circ\Delta(s)=\Phi(\iota_S(0,1))=\rho(0,1)=s
$$

$$
\Phi\circ\Delta(e_{[H]})=\Phi(\iota_D(e_{[H]}))=\pi(e_{[H]})=e_{[H]}
$$

and the other side
$$
\Delta\circ\Phi(\iota_S(g,n))=\Delta(\rho(g,n))=\Delta(u_gs^n)=\iota_S(g,n)
$$

$$
\Delta\circ\Phi(\iota_D(e_{[H]}))=\Delta(\pi(e_{[H]}))=\Delta(\iota_D(e_{[H]}))=\iota_D(e_{[H]}).
$$\par\vspace{\baselineskip}

Thus $\mathds{U}[\varphi,B]$ and $D[\varphi,B]\rtimes_{\alpha} S$ are isomorphic.
\begin{flushright}

  $\square$

  \end{flushright}

\begin{obs}\label{obs4}\em Note that $\mathds{U}[\varphi,B]$ is also isomorphic to $(D[\varphi,B]\rtimes_{\omega}G)\rtimes_{\tau}\mathds{N}$,
\begin{equation*}
\begin{split}
\omega: G &\rightarrow \hbox{Aut}(D[\varphi,B])\\
        g &\mapsto u_g(\cdot)u_{g^{-1}},\\ \\
\tau: \mathds{N} &\rightarrow \hbox{End}(D[\varphi,B]\rtimes_{\omega}G)\\
               n &\mapsto s^n(\cdot){s^*}^n
\end{split}
\end{equation*}

where for $a_g\delta_g$ of $D[\varphi,B]\rtimes_{\omega}G$, $\tau_n(a_g\delta_g)=s^na_g{s^*}^n\delta_{\varphi^n(g)}$.
\end{obs}

Using the minimal automorphic dilation presented by Laca in \cite{Laca1} it is possible to see the C$^*$-algebra $\mathds{U}[\varphi,B]$ as a corner of a group crossed product. For this, we need to prove the following.
\begin{prop}\label{propOre} The semidirect product $S=G\rtimes_\varphi\mathds{N}$ is an Ore semigroup i.e, it is cancellative and right-reversible.
\end{prop}
\dem Consider $(g_i,n_i)\in S$ for $i\in\{1,2,3\}$. $S$ is cancellative:
\begin{equation*}
\begin{split}
&(g_1,n_1)(g_3,n_3)=(g_2,n_2)(g_3,n_3)\\
\Rightarrow\;&(g_1\varphi^{n_1}(g_3),n_1+n_3)=(g_2\varphi^{n_2}(g_3),n_2+n_3)\\
\Rightarrow\; &n_1=n_2 \hbox{ and }g_1\varphi^{n_1}(g_3)=g_2\varphi^{n_1}(g_3)\\
\Rightarrow\; &g_1=g_2
\end{split}
\end{equation*}

\begin{equation*}
\begin{split}
&(g_1,n_1)(g_2,n_2)=(g_1,n_1)(g_3,n_3)\\
\Rightarrow\; &(g_1\varphi^{n_1}(g_2),n_1+n_2)=(g_1\varphi^{n_1}(g_3),n_1+n_3)\\
\Rightarrow\; &n_2=n_3 \hbox{ and }\varphi^{n_1}(g_2)=\varphi^{n_1}(g_3)\\
\Rightarrow\; &g_2=g_3\hbox{ as }\varphi\hbox{ is injective}.
\end{split}
\end{equation*}

Also any two principal left ideals of $S$ intersect:
\begin{equation*}
\begin{split}
(\varphi^{n_2}(g_1^{-1}),n_2)(g_1,n_1)&=(e,n_2+n_1)\\
&=(\varphi^{n_1}(g_2^{-1}),n_1)(g_2,n_2)\in S(g_1,n_1)\cap S(g_2,n_2).
\end{split}
\end{equation*}
\begin{flushright}

  $\square$

  \end{flushright}

It follows that the semigroup $S$ can be embedded in a group, called the enveloping group of $S$, which we will denote as $env(S)$, such that $S^{-1}S=env(S)$ (Theorem 1.1.2 \cite{Laca1}). It also implies that $S$ is a directed set by the relation defined by $(g,n)< (h,m)$ if $(h,m)\in S(g,n)$.  Let us define a candidate for $env(S)$. Consider
$$
\mathds{G}:=\D\lim_{\rightarrow}\{G_n:\varphi^n\}
$$

(with $G_n=G$ for all $n\in\mathds{N}$) and with the extended automorphism $\overline{\varphi}$ of $\mathds{G}$ construct the group
$$
\overline{S}:=\mathds{G}\rtimes_{\overline{\varphi}}\mathds{Z}.
$$

\begin{prop}\label{props}$\overline{S}\cong$ env$(S)$
\end{prop}
\dem For this we need to show that $S$ is a subsemigroup of $\overline{S}$ and $\overline{S}\subset S^{-1}S$ \cite{CliPre}.

First it is obvious that $S$ is a subsemigroup of the group $\overline{S}$ considering the inclusion $(g,n)\mapsto (g_0,n)$, where $g_0=g\in G=G_0\hookrightarrow\mathds{G}$.

Without loss of generality take $(g_i,j)\in\overline{S}$ with $i>|j|$. Then
$$
(g_i,j)=(g_i,-i)(e,j+i)=(g_0,i)^{-1}(e,j+i)\in S^{-1}S.
$$
\begin{flushright}

  $\square$

  \end{flushright}

Now consider the inductive system given by
$$
\overline{D}[\varphi,B]:=\D\lim_{\rightarrow}\left\{D[\varphi,B]_{(h,m)}:\alpha_{(h,m)}^{(g\varphi^n(h),n+m)}\right\}
$$

where
$$
D[\varphi,B]_{(h,m)}:=D[\varphi,B]
$$

and
$$
\alpha_{(h,m)}^{(g\varphi^n(h),n+m)}: D[\varphi,B]_{(h,m)}\rightarrow D[\varphi,B]_{(g,n)(h,m)}=D[\varphi,B]_{(g\varphi^n(h),n+m)},
$$\par\vspace{\baselineskip}

with $\alpha_{(h,m)}^{(g\varphi^n(h),n+m)}:=\alpha_{(g,n)}$ $\forall\,(h,m),(g,n)\in S$, where the latter was defined before Proposition \ref{prop3}. Then the C$^*$-dynamical system $(\overline{D}[\varphi,B],\overline{S},\overline{\alpha})$ is called the minimal automorphic dilation of $(D[\varphi,B],S,\alpha)$ where:
$$
\overline{\alpha}_{(g,n)}\circ\iota=\iota\circ\alpha_{(g,n)}, \forall\; (g,n)\in G\rtimes\mathds{N}
$$

with $\iota:D[\varphi,B]\hookrightarrow D[\varphi,B]_{(e,0)} \rightarrow \overline{D}[\varphi,B]$, and
$$
\overline{\D\bigcup_{(g,n)\in S}\overline{\alpha}_{(g,n)}^{-1}(\iota(D[\varphi,B]))}=\overline{D}[\varphi,B].
$$

Then by Theorem 2.2.1 in \cite{Laca1}:
\begin{lemma}\label{lema2} There exists an isomorphism
$$
\Phi:\mathds{U}[\varphi,B]\cong D[\varphi,B]\rtimes_{\alpha}S\cong\iota(1)(\overline{D}[\varphi,B]\rtimes_{\overline{\alpha}}\overline{S})\iota(1).
$$
\end{lemma}
\begin{flushright}

  $\square$

  \end{flushright}

Thus, $D[\varphi,B]\rtimes_{\alpha} S$ is Morita equivalent to $\overline{D}[\varphi,B]\rtimes_{\overline{\alpha}}\overline{S}$, $\Phi|_{D[\varphi,B]}=\iota$ and also $\Phi(u_gs^n)=\iota(1)\overline{U}_{(g,n)}\iota(1)$, where $\overline{U}:\overline{S}\rightarrow\mathcal{U}M(\overline{D}[\varphi,B]\rtimes_{\overline{\alpha}}\overline{S})$ (unitary multipliers).

\section{Separability, nuclearity and UCT}

By Proposition \ref{prop1} we conclude the following.
\begin{prop}\label{separab} If $B$ contains countably many subsets of $G$, then $\mathds{U}[\varphi,B]$ is separable.
\end{prop}
\dem With the condition satisfied we have countably many projections in $\mathds{U}[\varphi,B]$, and therefore it is generated by countably many elements.
\begin{flushright}

  $\square$

  \end{flushright}

And the group crossed product description obtained in last section implies two properties:
\begin{prop}\label{nuclear}If $G$ is amenable, $\mathds{U}[\varphi,B]$ is nuclear.
\end{prop}
\dem $G$ being amenable implies that $\overline{S}$ is amenable as well (amenability is closed under direct limits by \cite{vNeu} and also closed under semidirect products). But we know that $D[\varphi,B]$ is nuclear because it is commutative, therefore $\overline{D}[\varphi,B]\rtimes_{\overline{\alpha}}\overline{S}$ is nuclear by Proposition 2.1.2 in \cite{Ror}. Since hereditary C*-subalgebras of nuclear C*-algebras are nuclear by Corollary 3.3 (4) in \cite{ChoiEffros}, we conclude that
$$
\mathds{U}[\varphi,B]\cong D[\varphi,B]\rtimes_{\alpha}S\cong\iota(1)(\overline{D}[\varphi,B]\rtimes_{\overline{\alpha}}\overline{S})\iota(1)
$$

is nuclear.
\begin{flushright}

  $\square$

  \end{flushright}

\begin{prop}\label{uctprop}If $G$ is amenable, $\mathds{U}[\varphi,B]$ satisfies the UCT property.
\end{prop}
\dem Since $D[\varphi,B]$ is commutative, $\overline{D}[\varphi,B]\rtimes_{\overline{\alpha}}\overline{S}$ is isomorphic to a groupoid C$^*$-algebra. When the group $G$ is amenable then $\overline{S}$ also is, and the respective groupoid is also amenable. Therefore using a result by Tu (\cite{tutu} Proposition 10.7), the crossed product satisfies UCT. By Morita equivalence, $\mathds{U}[\varphi,B]$ also satisfies it.
\begin{flushright}

  $\square$

  \end{flushright}

\section{Purely infinite and simple}

To prove that under certain conditions our algebra is purely infinite and simple we use Proposition \ref{prop4} below, which is proven in Proposition 5.2 of \cite{Li1}.
\begin{prop}\label{prop4}Let $\widetilde{A}$ be a dense $*$-subalgebra of a unital C$^*$-algebra $A$. Assume that $\epsilon$ is a faithful conditional expectation on $A$ such that for every $0\neq x\in\widetilde{A}_{+}$ there exist finitely many projections $f_i\in A$ with
\begin{itemize}
  \item[(i)] $f_i\bot f_j$, $\,\forall\; i\neq j$,\vspace{0.2cm}
  \item[(ii)] $\exists\; s_i$ isometries such that $s_is_i^*=f_i$, $\,\forall\; i$,\vspace{0.2cm}
  \item[(iii)] $\left\|\D\sum_if_i\epsilon(x)f_i\right\|=\|\epsilon(x)\|$,\vspace{0.2cm}
  \item[(iv)] $f_ixf_i=f_i\epsilon(x)f_i\in\mathds{C}f_i$, $\,\forall\; i$.\vspace{0.2cm}
\end{itemize}
Then $A$ is purely infinite and simple.
\end{prop}
\begin{flushright}

  $\square$

  \end{flushright}

So we need to present a dense subalgebra and a conditional expectation of $\mathds{U}[\varphi,B]$. To find the conditional expectation, it is necessary to suppose in this section that the group $G$ is amenable. And to prove the main theorem, we suppose that $\varphi$ is pure, i.e:
$$
\D\bigcap_{n\in\mathds{N}}\varphi^n(G)=\{e\}.
$$

To start with, the next lemma tells us that
$$
S[\varphi,B]:=span(\{{s^*}^nu_{g^{-1}}e_{[H]}u_{g'}s^m:H\in B', g, g'\in G, n, m\in\mathds{N}\})
$$

is dense in $\mathds{U}[\varphi,B]$.
\begin{lemma}\label{lema3}The $*$-subalgebra of $\mathds{U}[\varphi,B]$ generated by
$$
\{e_{[H]}, u_g, s:\,H\in B, g\in G\}
$$

coincides with $S[\varphi,B]$.
\end{lemma}
\dem Note that
$$
\{e_{[H]},u_{g},s:\,H\in B, g\in G\}\subseteq S[\varphi,B]\subseteq span\{e_{[H]},u_g,s:H\in B, g\in G\}
$$

and $S[\varphi,B]$ is closed under multiplication:
\begin{equation*}
\begin{split}
&{s^*}^nu_{g^{-1}}e_{[H]}u_{g'}s^{n'}{s^*}^mu_{h^{-1}}e_{[K]}u_{h'}s^{m'}\\
&={s^*}^nu_{g^{-1}g'}(u_{g'^{-1}}e_{[H]}u_{g'})s^{n'}{s^*}^{n'}{s^*}^ms^{n'}(u_{h^{-1}}e_{[K]}u_h)u_{h^{-1}h'}s^{m'}\\
&={s^*}^nu_{g^{-1}g'}{s^*}^m(s^me_{[g'^{-1}H]}{s^*}^m)s^ms^{n'}{s^*}^{n'}{s^*}^m(s^{n'}e_{[h^{-1}K]}{s^*}^{n'})s^{n'}u_{h^{-1}h'}s^{m'}\\
&={s^*}^{n+m}u_{\varphi^m(g^{-1}g')}e_{[\varphi^m(g'^{-1}H)]}e_{[\varphi^{m+n'}(G)]}e_{[\varphi^{n'}(h^{-1}K)]}u_{\varphi^{n'}(h^{-1}h')}s^{n'+m'}\\
&={s^*}^{n+m}u_{\varphi^m(g^{-1}g')}e_{[\varphi^m(g'^{-1}H)\cap\varphi^{m+n'}(G)\cap\varphi^{n'}(h^{-1}K)]}u_{\varphi^{n'}(h^{-1}h')}s^{n'+m'}\\
&={s^*}^{n+m}u_{\varphi^m(g^{-1}g')}e_{[\varphi^m(g'^{-1})\varphi^m(H)\cap\varphi^{m+n'}(G)\cap\varphi^{n'}(h^{-1}K)]}u_{\varphi^{n'}(h^{-1}h')}s^{n'+m'}\\
&=\,0\in S[\varphi,B]\hbox{ or }\\
&={s^*}^{n+m}u_{\varphi^m(g^{-1}g')}e_{[\widetilde{g}(\varphi^m(H)\cap\varphi^{m+n'}(G)\cap\varphi^{n'}(K))]}
u_{\varphi^{n'}(h^{-1}h')}s^{n'+m'}\\
&={s^*}^{n+m}u_{\varphi^m(g^{-1}g')\widetilde{g}}e_{[\varphi^m(H)\cap\varphi^{m+n'}(G)\cap\varphi^{n'}(K))]}
u_{\widetilde{g}^{-1}\varphi^{n'}(h^{-1}h')}s^{n'+m'}\in S[\varphi,B].
\end{split}
\end{equation*}

The result follows.
\begin{flushright}

  $\square$

  \end{flushright}

Now we just have to define a conditional expectation to use in Proposition \ref{prop4} with the subalgebra defined above. For this we use the amenability of $G$. Therefore $\overline{S}$ is amenable, which implies that both the reduced and the full crossed products by $\overline{S}$ are isomorphic. With this we obtain a canonical conditional expectation of $\mathds{U}[\varphi,B]$.

Using the isomorphism
$$
\Phi:\mathds{U}[\varphi,B]\cong D[\varphi,B]\rtimes_{\alpha}S\rightarrow\iota(1)(\overline{D}[\varphi,B]\rtimes_{\overline{\alpha}}\overline{S})\iota(1).
$$

obtained in Lemma \ref{lema2} we have the easy-to-prove result below.
\begin{lemma}\label{lema4}There exists a faithful conditional expectation
\begin{equation*}
\begin{split}
\theta:\mathds{U}[\varphi,B]&\rightarrow\Phi^{-1}(i(1)\overline{D}[\varphi,B]i(1))\\
{s^*}^nu_{g^{-1}}e_{[H]}u_{g'}{s^{n'}}&\mapsto\left\{
                                                \begin{array}{ll}
                                                  {s^*}^nu_{g^{-1}}e_{[H]}u_gs^n, & \hbox{if }n=n'\hbox{ and }g=g'; \\
                                                  0, & \hbox{otherwise.}
                                                \end{array}
                                              \right.
\end{split}
\end{equation*}

for all $H\in B'$, $g,g'\in G$ and $n,n'\in\mathds{N}$.
\end{lemma}
\begin{flushright}

  $\square$

  \end{flushright}

Now we can prepare to prove that $\mathds{U}[\varphi,B]$ is simple and purely infinite, upon imposing some conditions. For this aim we follow and adapt the proof of Li \cite{Li1} (Section 5.2) and use the next lemmas to make the proof of the main theorem cleaner.

\begin{lemma}\label{lema5}Let $H$ and $G_i$ be distinct subgroups on $G$ with $\#\left[\dfrac{H}{H\cap G_i}\right]=\infty$ for all $1\leq i\leq n$. Then, for all $h, g_i\in G$, we have $hH\nsubseteq\D\bigcup_{i=1}^ng_i(H\cap G_i)$.
\end{lemma}
\dem By induction. For $n=1$:
$$
hH\subseteq g_1(H\cap G_1)\Rightarrow H\subseteq h^{-1}g_1(H\cap G_1)\Rightarrow \dfrac{H}{H\cap G_1}\neq\infty.
$$

Assume that the result holds for $n-1$. Let us prove it holds for $n$. Suppose that
$$
hH\subseteq\bigcup^n_{i=1}g_i(H\cap G_i),
$$

for some $h,g_i\in G$, with $1\leq i\leq n$. We can consider two possible cases:
\begin{itemize}
\item There exists $1<j\leq n$ with
$$
\#\left[\dfrac{H\cap G_1}{(H\cap G_1)\cap(H\cap G_j)}\right]<\infty.
$$
\end{itemize}

As
$$
\dfrac{(H\cap G_1)(H\cap G_j)}{H\cap G_j}\cong\dfrac{H\cap G_1}{(H\cap G_1)\cap(H\cap G_j)},
$$

it follows that the first one also has cardinality $<\infty$. But the exact sequence
$$
\dfrac{(H\cap G_1)(H\cap G_j)}{H\cap G_j}\hookrightarrow\dfrac{H}{H\cap G_j}\twoheadrightarrow\dfrac{H}{(H\cap G_1)(H\cap G_j)}
$$

with $\#\left[\frac{(H\cap G_1)(H\cap G_j)}{H\cap G_j}\right]<\infty$ and $\#\left[\frac{H}{H\cap G_j}\right]=\infty$ implies that
$$
\#\left[\dfrac{H}{(H\cap G_1)(H\cap G_j)}\right]=\infty.
$$

Define:
$$
\widetilde{G_i}:=\left\{
                  \begin{array}{ll}
                    H\cap G_i, & \hbox{if }G_i\neq G_1\hbox{ and }G_i\neq G_j; \\
                    (H\cap G_1)(H\cap G_j), & \hbox{if }G_i\in\{G_1,G_j\}.
                  \end{array}
                \right.
$$

Note that
$$
\#\left[\dfrac{H}{H\cap\widetilde{G_i}}\right]=\infty
$$

and
$$
hH\subseteq\bigcup^n_{i=1}g_i(H\cap G_i)\subseteq\bigcup^n_{i=1}g_i(H\cap\widetilde{G_i}),
$$

but the latter one contradicts our hypothesis, as $\#\{\widetilde{G_i}\}\leq n-1$.\par\vspace{\baselineskip}
\begin{itemize}
\item Now suppose that $\forall\; 1<j\leq n$,
$$
\#\left[\dfrac{H\cap G_1}{(H\cap G_1)(H\cap G_j)}\right]=\infty.
$$
\end{itemize}

As $\#[\frac{H}{H\cap G_1}]=\infty$, we have that $\exists\; g\in H$ such that $g(H\cap G_1)\neq g_i(H\cap G_i)$ $\forall\; 1\leq i\leq n$. Then:
\begin{equation*}
\begin{split}
g(H\cap G_1)&=g(H\cap G_1)\cap H\subseteq g(H\cap G_1)\cap\bigcup^n_{i=1}g_i(H\cap G_i)\\
&=\bigcup_{g(H\cap G_1)\cap g_i(H\cap G_i)\neq\emptyset}g(H\cap G_1)\cap g_i(H\cap G_i)\\
&=\bigcup_{g(H\cap G_1)\cap g_i(H\cap G_i)\neq\emptyset}\widetilde{g_i}((H\cap G_1)\cap(H\cap G_i))
\end{split}
\end{equation*}

and we can conclude that
$$
H\cap G_1\subseteq\bigcup_{g(H\cap G_1)\cap g_i(H\cap G_i)\neq\emptyset}g^{-1}\widetilde{g_i}((H\cap G_1)\cap(H\cap G_i)).
$$

But note that by construction $g(H\cap G_1)\cap g_i(H\cap G_1)=\emptyset$. So that union has been taken over less than $n$ elements, what contradicts our claim.
\begin{flushright}

  $\square$

  \end{flushright}

Let us show that $\mathds{U}[\varphi,B]$ together with the dense $*$-subalgebra $S[\varphi,B]$ (Lemma \ref{lema3}) and the faithful conditional expectation $\theta$ defined in Lemma \ref{lema4} satisfy the criteria of Proposition \ref{prop4}.

Take $0\neq x\in S[\varphi,B]_+$. As $\theta(x)\neq 0$, one has:
$$
\theta(x)=\D\sum^{{\scriptstyle{finite}}}_{(n',X)}\beta_{(n',X)}{s^*}^{n'}e_{[X]}s^{n'},
$$

where $(n',X)\in\mathds{N}\times C(B)$. Define $n$ to be the sum of all $n'$ with
$$
\beta_{(n',X)}{s^*}^{n'}e_{[X]}s^{n'}\neq 0.
$$

Then
$$
\theta(x)={s^*}^n\left(\D\sum^{{\scriptstyle{finite}}}_{(n',X)}\beta_{(n',X)}e_{[\varphi^{n-n'}(X)]}\right)s^n.
$$

Moreover using Lemma \ref{lema1}, it is possible to write
\begin{equation}\label{eq1}
\theta(x)={s^*}^n\left(\D\sum^{{\scriptstyle{finite}}}_{(g,H)}\beta_{(g,H)}e_{[gH]}\right)s^n,
\end{equation}

where the sum is over finitely many $(g,H)\in G\times B'.$\footnote{Remembering that $B'=\left\{\D\bigcap^m_{i=1}\varphi^{n_i}(H_i):H_i\in B\cup\{G\}, n_i\in\mathds{N}\right\}$.}

Note that
$$
{s^*}^ne_{[gH]}={s^*}^ns^n{s^*}^ne_{[gH]}={s^*}^ne_{[\varphi^n(G)]}e_{[gH]}={s^*}^ne_{[\varphi^n(G)\cap gH]},
$$

so we can assume that $gH\subseteq\varphi^n(G)$, for each $(g,H)\in G\times B'$.
\begin{lemma}\label{lema6} There exist finitely many pairwise orthogonal (nontrivial) projections $p_i$ in $\mathds{Z}$-$span(D[\varphi,B])$ such that $C^*(\{e_{[gH]}:\beta_{(g,H)}\neq 0\})=C^*(\{p_i\})$.
\end{lemma}
\dem Just orthogonalize the $e_{[g'H']}$. One can arrange the coefficients are integers.
\begin{flushright}

  $\square$

  \end{flushright}

Thus take some $p\in\{p_i\}$ among the $p_i$'s obtained above. Then
\begin{equation}\label{eq4}
p=\D\sum_jn_je_{[g_jH_j]}-\D\sum_{j'}\widetilde{n}_{j'}e_{[\widetilde{g}_{j'}\widetilde{H}_{j'}]}
\end{equation}

with finitely many $n_j,\widetilde{n}_{j'}\in\mathds{Z}_{>0}$ and $(g_j,H_j)$, $(\widetilde{g}_{j'},\widetilde{H}_{j'})\in G\times B'$.
\begin{lemma}\label{lema7}We can express $p$ as in (\ref{eq4}) so that $\forall\; K,\widetilde{K}\in\{H_j,\widetilde{H}_j\}$ the cardinality of
$\frac{K}{K\cap\widetilde{K}}$ is 1 or $\infty$.
\end{lemma}
\dem By induction. Enumerate $\{H_j,\widetilde{H}_j\}$ by $\{K_i\}$. Of course the lemma holds if there is just $K_1$.

Suppose that it holds for $\{K_1,\ldots,K_h\}$. Define $K_{h+1}^{(0)}:=K_{h+1}$ and for\nl $j=1,\dots h$
\begin{equation}\label{eq5}
K_{h+1}^{(j)}:=\left\{
                \begin{array}{ll}
                  K_{h+1}^{(j-1)}, & \hbox{if }\#[K_{h+1}^{(j-1)}/(K_{h+1}^{(j-1)}\cap K_j)]\in\{1,\infty\}, \\
                  K_{h+1}^{(j-1)}\cap K_j, & \hbox{otherwise.}
                \end{array}
              \right.
\end{equation}

We want to change $K_{h+1}$ successively to $K_{h+1}^{(0)}, K_{h+1}^{(1)}, \cdots$, until $K_{h+1}^{(h)}$.

Suppose that $K_{h+1}^{(j)}= K_{h+1}^{(j-1)}\cap K_j$ as described above in (\ref{eq5}). Therefore we have $1<\#[K_{h+1}^{(j-1)}/(K_{h+1}^{(j-1)}\cap K_j)]=M<\infty$, and then, $K_{h+1}^{(j-1)}=\D\bigcup_{i=1}^M g_i(K_{h+1}^{(j-1)}\cap K_j)$.

So we can replace $K_{h+1}$ by $K_{h+1}':=K^{(h)}_{h+1}$, because the projections will still be written using the initial $\{K_i\}$.\par\vspace{\baselineskip}

\afi
$$\#\left[\dfrac{K_{h+1}'}{K_{h+1}'\cap K}\right]\in\{1,\infty\}, \forall\; K\in\{K_1,\cdots,K_h\}.
$$

\prafi Let us prove by induction on $j$ that $\#\left[\dfrac{K_{h+1}^{(j)}}{K_{h+1}^{(j)}\cap K}\right]\in\{1,\infty\}$,\par\vspace{\baselineskip}\par\vspace{\baselineskip} for every $K\in\{K_1,\cdots,K_j\}$. By construction it holds for $j=1$. Suppose it holds for $j-1$, that is
$$
\#\left[\dfrac{K_{h+1}^{(j-1)}}{K_{h+1}^{(j-1)}\cap K}\right]\in\{1,\infty\},\,\forall\; K\in\{K_1,\cdots,K_{j-1}\},
$$

and let us prove the assertion for $j$. Also by construction $\#\left[\dfrac{K_{h+1}^{(j)}}{K_{h+1}^{(j)}\cap K_j}\right]$ belongs to $\{1,\infty\}$. Then, we need to show that
$$
\#\left[\dfrac{K_{h+1}^{(j)}}{K_{h+1}^{(j)}\cap K}\right]\in\{1,\infty\},\;\forall\; K\in\{K_1,\cdots,K_{j-1}\}.
$$

If $K_{h+1}^{(j)}=K_{h+1}^{(j-1)}$, then this holds by the induction hypothesis.

But, if $K_{h+1}^{(j)}=K_{h+1}^{(j-1)}\cap K_j$, then $K_{h+1}^{(j)}\subset K_{h+1}^{(j-1)}$ and therefore it follows that 
$$
1<\#\left[\dfrac{K_{h+1}^{(j-1)}}{K_{h+1}^{(j-1)}\cap K_j}\right]<\infty.
$$

Now, by our induction hypothesis, we have two possibilities for each\\
$K\in\{K_1,\ldots,K_{j-1}\}$:

$\;\;\;\;\;\;\;\;\bullet\;\#\left[\dfrac{K_{h+1}^{(j-1)}}{K_{h+1}^{(j-1)}\cap K}\right]=1$: in this case, as $K_{h+1}^{(j)}\subset K_{h+1}^{(j-1)}\subset K$, it follows \par\vspace{\baselineskip} that $\#\left[\dfrac{K_{h+1}^{(j)}}{K_{h+1}^{(j)}\cap K}\right]=1$.\par\vspace{\baselineskip}

$\;\;\;\;\;\;\;\;\bullet\;\#\left[\dfrac{K_{h+1}^{(j-1)}}{K_{h+1}^{(j-1)}\cap K}\right]=\infty$.\par\vspace{\baselineskip}

Consider the exact sequence:
$$
\dfrac{K_{h+1}^{(j)}}{K_{h+1}^{(j)}\cap K}\hookrightarrow\dfrac{K_{h+1}^{(j-1)}}{K_{h+1}^{(j)}\cap K}\twoheadrightarrow\dfrac{K_{h+1}^{(j-1)}}{K_{h+1}^{(j)}}.
$$

The inclusion $\dfrac{K_{h+1}^{(j-1)}}{K_{h+1}^{(j-1)}\cap K}\subset\dfrac{K_{h+1}^{(j-1)}}{K_{h+1}^{(j)}\cap K}$ implies that the second term has size $\infty$. The third term has cardinality $<\infty$ because it is equal to $\dfrac{K_{h+1}^{(j-1)}}{K_{h+1}^{(j-1)}\cap K_j}$. As that sequence is exact, we must have $\#\left[\dfrac{K_{h+1}^{(j)}}{K_{h+1}^{(j)}\cap K}\right]=\infty$. \par\vspace{\baselineskip}

Thus we conclude that $\#\left[\dfrac{K_{h+1}^{(j)}}{K_{h+1}^{(j)}\cap K}\right]\in\{1,\infty\}$, $\forall\; K\in\{K_1,\cdots,K_j\}$.
\begin{flushright}

  $\square$

  \end{flushright}

Set
\begin{equation}\label{eq7}
K_j':=\left\{
                \begin{array}{ll}
                  K_j\cap K_{h+1}', & \hbox{if }1<\#\left[\dfrac{K_j}{K_j\cap K_{h+1}'}\right]<\infty, \\
                  K_j, & \hbox{otherwise}
                \end{array}
              \right.
\end{equation}

for $j=1,\dots h$. This gives a new sequence $\{K_1',\cdots,K_{h+1}'\}$. And then it only remains to prove that
\begin{equation*}
\#\left[\dfrac{K_j'}{K_j'\cap K_{\widetilde{j}}'}\right]\in\{1,\infty\}.
\end{equation*}

Note that, if $j$ or $\widetilde{j}$ is equal to $h+1$, this holds by the claim above.

So, suppose that $j$ and $\widetilde{j}$ are in $\{1,\cdots,h\}$. Then, by our induction hypothesis, we have two possibilities:\par\vspace{\baselineskip}

$\;\;\;\;\;\;\;\;\bullet$ $\#\left[\dfrac{K_j}{K_j\cap K_{\widetilde{j}}}\right]=1$:
then $K_j\subseteq K_{\widetilde{j}}$.\par\vspace{\baselineskip}

If $K_j'=K_j\cap K_{h+1}'$, then $K_j'\subseteq K'_{\widetilde{j}}$, and (\ref{eq7}) holds.

Otherwise, $K_j'=K_j\neq K_j\cap K_{h+1}'$, and therefore $\#\left[\dfrac{K_j}{K_j\cap K_{h+1}'}\right]=\infty$. Then (as $K_j\subseteq K_{\widetilde{j}}$) we have the inclusion:
$$
\dfrac{K_j}{K_j\cap K_{h+1}'}\subseteq\dfrac{K_{\widetilde{j}}}{K_{\widetilde{j}}\cap K_{h+1}'},
$$

which implies $\left[\dfrac{K_{\widetilde{j}}}{K_{\widetilde{j}}\cap K_{h+1}'}\right]=\infty$. So $K'_{\widetilde{j}}=K_{\widetilde{j}}$ and our claim holds.

$\;\;\;\;\;\;\;\;\bullet$ $\#\left[\dfrac{K_j}{K_j\cap K_{\widetilde{j}}}\right]=\infty\;$:
As $\dfrac{K_j}{K_j\cap K_{\widetilde{j}}}\subseteq\dfrac{K_j}{K_j\cap K_{\widetilde{j}}'}$, if $K_j'=K_j$ the claim holds.

Now, if $K_j'=K_j\cap K_{h+1}'\neq K_j$, then we have the exact sequence:
$$
\dfrac{K_j'}{K_j'\cap K_{\widetilde{j}}}\hookrightarrow\dfrac{K_j}{K_j'\cap K_{\widetilde{j}}}\twoheadrightarrow\dfrac{K_j}{K_j'}.
$$

The set $\dfrac{K_j}{K_j\cap K_{\widetilde{j}}}$ has size $\infty$ and is contained in the second term, so it has size $\infty$ too. The third term has size $<\infty$ as $K_j\cap K_{h+1}'\neq K_j$ implies that $\#\left[\dfrac{K_j}{K_j\cap K_{h+1}'}\right]<\infty$.

Hence, we conclude that $\left[\dfrac{K_j'}{K_j'\cap K_{\widetilde{j}}}\right]=\infty$, proving the lemma.
\begin{flushright}

  $\square$

  \end{flushright}

\begin{lemma}\label{lema8} There exist finitely many pairwise orthogonal projections\nl $p_i\in\mathds{U}[\varphi,B]$ such that
$$
C^*(\{p_i\})\cong C^*(\{e_{[gH]}:\beta_{(g,H)}\neq 0\}),
$$

where the $(g,H)$'s come from equation (\ref{eq1}). Moreover if exists $m\in\mathds{N}$ such that $\varphi^m(G)\subseteq\bigcap_{H\in B}H$ then for all $i$, there exists $h_i\in G$ and $m_i\in\mathds{N}$ with
$$
e_{[h_i\varphi^{m_i}(G)]}\leq p_i.
$$
\end{lemma}
\dem We have
$$
\theta(x)={s^*}^n\left(\D\sum_{(g,H)}^{\hbox{\tiny finite}}\beta_{(g,H)}e_{[gH]}\right)s^n,\hbox{ with }(g,H)\in G\times B',
$$

where we recall that
$$
B'=\left\{\D\bigcap_{i=1}^m\varphi^{n_i}(H_i):H_i\in B\cup\{G\}, n_i\in\mathds{N}\right\}.
$$

We can assume that $gH\subset\varphi^n(G)$ and, by Lemma \ref{lema6} we have finitely many pairwise orthogonal projections $p_i$ in $\mathds{Z}$-$span(D[\varphi,B])$ with
$$
C^*(\{e_{[gH]}:\beta_{(g,H)}\neq 0\})=C^*(\{p_i\}).
$$

Choose some $p\in\{p_i\}$ and write it as
$$
p=\D\sum_jn_je_{[g_jH_j]}-\D\sum_{j'}\widetilde{n}_{j'}e_{[\widetilde{g}_{j'}\widetilde{H}_{j'}]}
$$

with finitely many $n_j, \widetilde{n}_{j'}\in\mathds{Z}_{>0}$. We can write $p$ such that each projection $e_{[g,H]}$ appears at most one time and $\#\left[\dfrac{K}{K\cap \widetilde{K}}\right]\in\{1,\infty\}$  $\forall\; K, \widetilde{K}\in\{H_j,\widetilde{H}_{j'}\}$ by Lemma \ref{lema7}.

Choose some maximal $H\in\{H_j,\widetilde{H}_{j'}\}$.

Take $g\in G$ and $n\in\mathds{Z}_{>0}$ so that $ne_{[gH]}$ appears in $p$. Multiplying $p$ with $e_{[gH]}$ gives
\begin{equation*}
e_{[gH]}p=ne_{[gH]}+\D\sum_kn_ke_{[c_k(H\cap H_k)]}-\D\sum_l\widetilde{n}_le_{[\widetilde{c}_l(H\cap\widetilde{H}_l)]},
\end{equation*}

for (finitely many) $c_k,\widetilde{c}_l\in G$ and $n_k, \widetilde{n}_l\in\mathds{Z}_{>\,0}$.

Note that we must have $\#\left[\dfrac{H}{H\cap H_k}\right]=\infty$ because if $\#\left[\dfrac{H}{H\cap H_k}\right]=1$ then $H_k= H$ would imply $e_{[g_jH_j]}= e_{[\widetilde{g}_{j'}\widetilde{H}_{j'}]}$ for some $j$ and $j\,'$.

Then, by Lemma \ref{lema5},
$$
gH\nsubseteq\left[\bigcup_kc_k(H\cap H_k)\right]\cup\left[\bigcup_l\widetilde{c}_l(H\cap \widetilde{H}_l)\right],
$$

which allows us to find $r\in gH\backslash\left[\D\bigcup_kc_k(H\cap H_k)\right]\cup\left[\D\bigcup_l\widetilde{c}_l(H\cap \widetilde{H}_l)\right]$.

One can conclude that:
\begin{equation*}
\begin{split}
e_{[r(\cap_k(H\cap H_k)\cap\cap_l(H\cap\widetilde{H}_l))]}&\leq e_{[gH]},\\
e_{[r(\cap_k(H\cap H_k)\cap\cap_l(H\cap\widetilde{H}_l))]}&\perp e_{[c_k(H\cap H_k)]},\;\forall\; k, and\\
e_{[r(\cap_k(H\cap H_k)\cap\cap_l(H\cap\widetilde{H}_l))]}&\perp e_{[\widetilde{c}_l(H\cap\widetilde{H}_l)]},\;\forall\; l.
\end{split}
\end{equation*}

Multiplying the equation above by $e_{[r(\cap_k(H\cap H_k)\cap\cap_l(H\cap\widetilde{H}_l))]}$ gives
\begin{equation*}
\begin{split}
e_{[r(\cap_k(H\cap H_k)\cap\cap_l(H\cap\widetilde{H}_l))]}p=ne_{[r(\cap_k(H\cap H_k)\cap\cap_l(H\cap\widetilde{H}_l))]}.
\end{split}
\end{equation*}

As the first term is a projection (because it is the product of two commuting projections) we must have $n=1$. So, $e_{[r(\cap_k(H\cap H_k)\cap\cap_l(H\cap\widetilde{H}_l))]}\leq p$.

If our additional hypothesis is satisfied, we have $\widetilde{m}\in\mathds{N}$ such that\footnote{$\varphi^m(G)\subseteq\displaystyle\bigcap_{H\in B}H\Rightarrow \varphi^{\widetilde{m}}(G)\subseteq\bigcap_{H_i\in B'\atop 0\leq i\leq n}H_i$, for some $\widetilde{m}$ bigger than $m$.}
$$
e_{[r\varphi^{\widetilde{m}}(G)]}\leq e_{[r(\cap_k(H\cap H_k)\cap\cap_l(H\cap\widetilde{H}_l))]}\leq p.
$$

therefore we just have to denote $h_i=r$ and $m_i=\widetilde{m}$. The conclusion holds if this is done for every element of $\{p_i\}$.
\begin{flushright}

  $\square$

  \end{flushright}

\begin{obs}\em\label{obs5}  Note that in last lemma for every $i$ we can choose $m_i$ as big as we want, because $\varphi^{m+1}(G)\subset\varphi^m(G)$.
\end{obs}

\begin{teo}\label{teo1}Let $G$ be an amenable group, $B$ some family of subgroups
in $G$ containing $G$ and $\varphi$ a pure injective endomorphism of $G$. Also suppose that $\exists\; k\in\mathds{N}$ such that $\varphi^k(G)\subseteq\bigcap_{H\in B}H$.

Then the C$^*$-algebra $\mathds{U}[\varphi,B]$ is purely infinite and simple.
\end{teo}
\dem We already have the candidates to use with Proposition \ref{prop4}, namely
$$
\theta:\mathds{U}[\varphi,B]\rightarrow\Phi^{-1}(\iota(1)\overline{D}[\varphi,B]\iota(1)),
$$

and
$$
S[G,B]=span(\{{s^*}^nu_{g^{-1}}e_{[I]}u_{g'}s^m:I\in B, g,g'\in G,n,m\in\mathds{N}\}).
$$

Take $0\neq x\in S[G,B]_{sa}$. Then
$$
x=\D\sum_{(g,g',l,l',J)}\alpha_{(g,g',l,l',J)}{s^*}^lu_{g^{-1}}e_{[J]}u_{g'}s^{l'}.
$$

As in previous Lemma \ref{lema8},
$$
\theta(x)={s^*}^n\left(\D\sum_{(g,H)}^{\hbox{{\tiny finite}}}\beta_{(g,H)}e_{[gH]}\right)s^n,
$$

for some $n\in\mathds{N}$ and $(g,H)\in G\times B'$ with $\beta_{(g,H)}\neq 0$ where $gH\subset\varphi^n(G)$.

By Lemma \ref{lema8} we find finitely many pairwise orthogonal (nontrivial) projections $\{p_i\}$ with $C^*(\{e_{[gH]}:\beta_{(g,H)}\neq 0\})=C^*(\{p_i\})$ and there exist $m_i\in\mathds{N}$ and $h_i\in G$ such that $e_{[h_i\varphi^{m_i}(G)]}\leq p_i\leq e_{[\varphi^n(G)]}$ $\forall\; i$. Using Remark \ref{obs5}, we can suppose that $m_i\geq n$ $\forall\; i$. Also note that $h_i\in\varphi^n(G)$.

Thus the projections $F_i:={s^*}^ne_{[h_i\varphi^{m_i}(G)]}s^n$ satisfy $F_i\leq {s^*}^np_i s^n$ and
\begin{equation}\label{equa113}
\begin{split}
C^*(\{{s^*}^ne_{[gH]}s^n:\beta_{(g,H)}\neq 0\})=C^*(\{{s^*}^np_is^n\})&\rightarrow C^*(\{F_i\})\\
y&\mapsto\D\sum_iF_iyF_i
\end{split}
\end{equation}

is an isomorphism that maps ${s^*}^np_is^n$ to $F_i$.

These projections $F_i$ satisfy only (i) and (ii) of the conditions in Proposition \ref{prop4}.

Call $(g,g',l,l',J)$ critical if $\alpha_{(g,g',l,l',J)}{s^*}^lu_{g^{-1}}e_{[J]}u_{g'}s^{l'}\neq 0$ and $\delta_{g,g'}\delta_{l,l'}=0$. Note that
$$
x-\theta(x)=\sum_{(g,g',l,l',J)\hbox{ {\tiny critical}}}{s^l}^*u_{g^{-1}}e_{[J]}u_{g'}s^{l'}.
$$

But for each $i$, it is possible to take some $a_i\in\varphi^{-n}(h_i)\varphi^{m_i-n}(G)$ 
satisfying 
$$
\varphi^{l'}(a_i^{-1})g'^{-1}g\varphi^l(a_i)\neq e
$$ 
for all critical $(g,g',l,l',J)$.

Surely, if not then we have $r_1\neq r_2\in\varphi^{m_i-n}(G)$ such that
$$
\varphi^{l'}(r_1^{-1})g'^{-1}g\varphi^l(r_1)=e=\varphi^{l'}(r_2^{-1})g'^{-1}g\varphi^l(r_2).
$$

If $l=l'$ we have $g\neq g'$ (as $\delta_{g,g'}\delta_{l,l'}=0$) and then
$$
\varphi^{l}(r_1^{-1})g'^{-1}g\varphi^l(r_1)=e\Rightarrow g'^{-1}g=e
$$

which contradicts $g\neq g'$.

Suppose now that $l\neq l'$. As $r_1=r_2r_2^{-1}r_1$ we get
$$
e=\varphi^{l'}((r_2r_2^{-1}r_1)^{-1})g'^{-1}g\varphi^l(r_2r_2^{-1}r_1)=\varphi^{l'}(r_1^{-1}r_2)\varphi^l(r_2^{-1}r_1)
$$

which implies that $r_1=r_2$ (because $\varphi$ is pure). This contradicts our assumptions.

Now as our endomorphism $\varphi$ is pure, for all critical $(g,g',l,l',J)$ and for all $i$ there exists $n_{(g,g',l,l',J,i)}\in\mathds{N}$ (as big as we need) such that $\varphi^{l'}(a_i^{-1})g'^{-1}g\varphi^l(a_i)\notin\varphi^{n_{(g,g',l,l',J,i)}}(G)$. Let us call
\begin{equation*}
b_i:=(m_i-n)\prod_{(g,g',l,l',J)\hbox{ {\tiny critical}}}n_{(g,g',l,l',J,i)}.
\end{equation*}

Note that
\begin{equation}\label{eq10}
\varphi^{l'}(a_i^{-1})g'^{-1}g\varphi^l(a_i)\notin\varphi^{b_i}(G).
\end{equation}

Define $f_i:=e_{[a_i\varphi^{b_i}(G)]}$. We want to prove that these projections satisfy the conditions of Proposition \ref{prop4}, which are:
\begin{itemize}
  \item[(i)] $f_i\bot f_j$,$\forall\; i\neq j$,\par\vspace{\baselineskip}
  \item[(ii)] $f_i\sim_{z_i} 1$, via isometries $z_i\in A$,$\forall\; i$,\par\vspace{\baselineskip}
  \item[(iii)] $\left\|\D\sum_if_i\theta(x)f_i\right\|=\|\theta(x)\|$, and\par\vspace{\baselineskip}
  \item[(iv)] $f_ixf_i=f_i\theta(x)f_i\in\mathds{C}f_i$,$\forall\; i$.
\end{itemize}

As $b_i\geq m_i-n$ and $\varphi^n(a_i)\in h_i\varphi^{m_i}(G)$ it follows that $s^ne_{[a_i\varphi^{b_i}(G)]}{s^*}^n\leq e_{[h_i\varphi^{m_i}(G)]}$ and then
$$
f_i={s^*}^ns^ne_{[a_i\varphi^{b_i}(G)]}{s^*}^ns^n\leq {s^*}^ne_{[h_i\varphi^{m_i}(G)]}s^n=F_i.
$$

This implies that $f_i\perp f_j$ $\forall\; i\neq j$ are pairwise orthogonal and (i) is satisfied. Item (ii) is also easily satisfied, because
$$
f_i=e_{[a_i\varphi^{b_i}(G)]}=(u_{a_i}s^{b_i})(u_{a_i}s^{b_i})^*\sim (u_{a_i}s^{b_i})^*(u_{a_i}s^{b_i})=1.
$$

As (\ref{equa113}) is an isomorphism and $f_i\leq F_i$ the map
\begin{equation*}
\begin{split}
C^*(\{{s^*}^ne_{[gH]}s^n:\beta_{(g,H)}\neq 0\})&\rightarrow C^*(\{f_i\})\\
y&\mapsto\D\sum_if_iyf_i
\end{split}
\end{equation*}

is an isomorphism as well. Therefore it is isometric and (iii) is satisfied.

And finally, for the last condition, let us expand $f_i(x-\theta(x))f_i$:
\begin{equation*}
\begin{split}
&f_i(x-\theta(x))f_i\\
&=f_i\left(\D\sum_{(g,g',l,l',J)\hbox{ {\tiny critical}}}\beta_{(g,g',l,l',J)}{s^*}^lu_{g^{-1}}e_{[J]}u_{g'}s^{l'}\right)f_i\\
&=\D\sum_{(g,g',l,l',J)\hbox{{\tiny critical}}}\beta_{(g,g',l,l',J)}{s^*}^lu_{g^{-1}}(u_{g}s^{l}f{s^*}^lu_{g^{-1}})e_{[J]}(u_{g'}s^{l'}f{s^*}^{l'}u_{{g'}^{-1}})u_{g'}s^{l'}\\
&=\D\sum_{(g,g',l,l',J)\hbox{{\tiny critical}}}\beta_{(g,g',l,l',J)}{s^*}^lu_{g^{-1}}e_{[g\varphi^l(a_i)\varphi^{l+b_i}(G)]}e_{[g'\varphi^{l'}(a_i)\varphi^{l'+b_i}(G)]}
e_{[J]}u_{g'}s^{l'}.
\end{split}
\end{equation*}

Now, note that
$$
[g\varphi^l(a_i)\varphi^{l+b_i}(G)]\cap[g'\varphi^{l'}(a_i)\varphi^{l'+b_i}(G)]
\neq\emptyset\Rightarrow\varphi^{l'}(a_i^{-1}){g'}^{-1}g\varphi^l(a_i)\in\varphi^{b_i}(G)
$$

which is a contradiction with our choice of $b_i$ by (\ref{eq10}). So the intersection above must be empty and then $f_ixf_i=f_i\theta(x)f_i\in\mathds{C}f_i$, $\forall\; i$.

Therefore, by Proposition \ref{prop4}, our C$^*$-algebra is simple and purely infinite.
\begin{flushright}

  $\square$

  \end{flushright}

\begin{corol} When satisfied the conditions of the theorem above, the concrete C$^*$-algebra $C_r^*[\varphi,B]$ is isomorphic to the universal one $\mathds{U}[\varphi,B]$, as defined in Definitions \ref{defi1} and \ref{defi2} respectively.
\end{corol}
\begin{flushright}

  $\square$

  \end{flushright}

\begin{teo}\label{teokirch2} If the conditions of the theorem above are satisfied, the universal C$^*$-algebra $\mathds{U}[\varphi,B]$ is a Kirchberg algebra satisfying the UCT property (Propositions \ref{separab}, \ref{nuclear} and \ref{uctprop}).
\end{teo}
\begin{flushright}

  $\square$

  \end{flushright}

\section{The case $B=\{G\}$}\label{seccan}

In the following chapter we study the particular case when $B$ contains only subgroups of the form $g\varphi^k(G)$ for 
$k\in\mathds{N}$ and $g\in G$. It will be now denoted $\mathds{U}[\varphi]$ and its K-theory will be calculated using a similar idea as presented in \cite{CunVer} i.e, using the continuity of the functors $K_0$ and $K_1$, and also the Khoshkam-Skandalis sequence \cite{Khoska}. We conclude that $K_*(\mathds{U}[\varphi])\cong K_*(C^*(G))$ and that when $G$ is amenable, the C$^*$-algebras $\mathds{U}[\varphi]$ are classifiable by Kirchberg's classification theorem.

To finish we use the recently-introduced semigroup C$^*$-algebras from \cite{Li0} and \cite{Li2} and show that $\mathds{U}[\varphi]$ is isomorphic to the full semigroup C$^*$-algebra of the semigroup $S=G\rtimes_\varphi\mathds{N}$. This implies that when the group $G$ is amenable and the endomorphism $\varphi$ is pure the three semigroup C$^*$-algebras (the full one, the reduced one and a third one, given by viewing $S$ as a subsemigroup of the group $\overline{S}$ defined before Proposition \ref{props}) defined by Li are isomorphic to $\mathds{U}[\varphi]$ and classifiable by Kirchberg's classification theorem.

By Proposition \ref{prop2} if we choose $\overline{B}$ containing subsets of the form $g\varphi^k(G)$ for $k\in\mathds{N}$ and $g\in G$, then $\mathds{U}[\varphi,B]$ is isomorphic to the one obtained when we start only with $B=\{G\}$. Therefore:
\begin{prop} \label{prop5}When $B$ contains only subsets of the form $g\varphi^k(G)$, for some fixed $\varphi$, $k\in\mathds{N}$ and $g\in G$, the C$^*$-algebra $\mathds{U}[\varphi,B]$ is isomorphic to $\mathds{U}[\varphi]$ and can be redefined as the universal C$^*$-algebra generated by
\begin{equation*}
\begin{split}
\hbox{unitaries }&\{u_g:g\in G\}\\
\hbox{and one isometry }&\{s\}
\end{split}
\end{equation*}
satisfying:
\begin{enumerate}
	\item[(i)] $u_gs^nu_hs^m=u_{g\varphi^n(h)}s^{n+m}$;\par\vspace{\baselineskip}
    \item[(ii)] $u_gs^n{s^*}^nu_{g^{-1}}u_hs^m{s^*}^mu_{h^{-1}}=u_hs^m{s^*}^mu_{h^{-1}}u_gs^n{s^*}^nu_{g^{-1}}=\nl \left\{
    \begin{array}{ll}
      u_gs^n{s^*}^nu_{g^{-1}}, & \hbox{if }h\in g\varphi^m(G); \\
      0, & \hbox{otherwise,}
    \end{array}
  \right.$\par\vspace{\baselineskip}
for $n\geq m$.
\end{enumerate}
\end{prop}
\begin{flushright}

  $\square$

  \end{flushright}

A simple use of Propositions \ref{separab}, \ref{nuclear} and \ref{uctprop}, and Theorem \ref{teo1} gives the following.

\begin{prop} The C$^*$-algebra $\mathds{U}[\varphi]$ is separable. When the group $G$ is amenable, it is also nuclear and satisfies UCT. Furthermore if $G$ is amenable and $\varphi$ is pure, then $\mathds{U}[\varphi]$ is also simple and purely infinite, therefore a Kirchberg algebra satisfying UCT.
\end{prop}
\begin{flushright}

  $\square$

  \end{flushright}

\subsection{K-theory}

Using Remark \ref{obs4} we know that
$$
\mathds{U}[\varphi]\cong (D[\varphi]\rtimes_{\omega}G)\rtimes_{\tau}\mathds{N}
$$

with
\begin{equation*}
\begin{split}
\omega: G &\rightarrow \hbox{Aut}(D[\varphi])\\
        g &\mapsto u_g(\cdot)u_{g^{-1}},\\ \\
\tau: \mathds{N} &\rightarrow \hbox{End}(D[\varphi]\rtimes_{\omega}G)\\
               n &\mapsto s^n(\cdot){s^*}^n
\end{split}
\end{equation*}

where for $a_g\delta_g\in D[\varphi]\rtimes_{\omega}G$, $\tau_n(a_g\delta_g)=s^na_g{s^*}^n\delta_{\varphi^n(g)}$. But note that
$$
D[\varphi]\cong\D\lim_{\rightarrow \atop n}D_n
$$

for $n\in\mathds{N}$ with
$$
D_n:=C^*\left(\left\{u_gs^k{s^*}^ku_{g^{-1}}:\; 0\leq k\leq n,\; g\in\frac{G}{\varphi^k(G)}\right\}\right)
$$

and the inclusion being the identity. Therefore
$$
D[\varphi]\rtimes_{\omega}G\cong\D\lim_{\rightarrow\atop n}(D_n\rtimes_{\omega}G),
$$

where
$$
D_n\rtimes_{\omega}G\cong C^*\left(\left\{u_gs^k{s^*}^ku_{h^{-1}}:\; 0\leq k\leq n,\; g,h\in G\right\}\right).
$$

Moreover, for $k\in\mathds{N}$, $A_k:=C^*(\{u_gs^k{s^*}^ku_{h^{-1}}:\,g,h\in G\})$ is an ideal of $D_k\rtimes_{\omega}G$, because for $m\leq k$,
$$
u_hs^k{s^*}^ku_{h^{-1}}u_gs^m{s^*}^mu_{g^{-1}}=\left\{
                                                 \begin{array}{ll}
                                                   u_hs^k{s^*}^ku_{h^{-1}}, & \hbox{if }g\in h\varphi^m(G)\\
                                                                         0, & \hbox{otherwise.}
                                                 \end{array}
                                               \right.
$$

But note that every element $u_gs^k{s^*}^ku_{{h}^{-1}}$ in $A_k$ can be uniquely written as\\
$u_{g_i}s^k{s^*}^ku_{{g_j}^{-1}}u_{\varphi^k(t)}$, for $g_i,g_j\in\frac{G}{\varphi^k(G)}$ and $t\in G$. Therefore, if one defines the correspondence
$$
u_gs^k{s^*}^ku_{{h}^{-1}}=u_{g_i}s^k{s^*}^ku_{{g_j}^{-1}}u_{\varphi^k(t)}\mapsto E_{i,j}\otimes u_{\varphi^k(t)},
$$

where $\{E_{i,j}\}$ is the family of unit matrices which give rise to the set $\mathbf{K}$ of compact operators, it follows that
$$
A_k\cong\mathbf{K}\otimes C^*(\varphi^k(G))\cong\mathbf{K}\otimes C^*(G).
$$

So starting with the case $n=1$, we can build the following exact sequence:
$$
0\rightarrow A_1\xrightarrow{\iota} D_1\rtimes_{\omega}G\xrightarrow{\rho} C^*(G)\rightarrow 0
$$

where ${\iota}$ and ${\rho}$ are the canonical inclusion and projection maps respectively. But the sequence above splits if we also consider the canonical inclusion
$$
\gamma: C^*(G)\rightarrow D_1\rtimes_{\omega}G.
$$

This implies that the corresponding exact sequence of K-groups also splits, which means that (using the K\"{u}nneth Formula \cite{Scho})
\begin{equation*}
\begin{split}
K_{*}(D_1\rtimes_{\omega}G)&\cong K_{*}(C^*(G))\oplus K_{*}(\mathbf{K}\otimes C^*(\varphi(G)))\\
&\cong K_{*}(C^*(G))\oplus K_{*}(C^*(G)).
\end{split}
\end{equation*}

Using the same argument repeatedly, it is easy to conclude that
$$
K_{*}(D_n\rtimes_{\omega}G)=\bigoplus_{i=0}^n K_{*}(C^*(G))
$$

and consequently
\begin{equation}\label{eq11}
K_{*}(D[\varphi]\rtimes_{\omega}G)=\lim_{\rightarrow \atop n}\bigoplus_{i=0}^n K_{*}(C^*(G))=\bigoplus_{\mathds{N}} K_{*}(C^*(G)),
\end{equation}

where the $k$-th group $K_{*}(C^*(G))$ of the direct sum above represents the K-group of $A_k$.

Applying the Khoshkam-Skandalis sequence for $\mathds{N}$-crossed products \cite{Khoska}, we have the following sequence:
$$
\begin{array}{ccccc}
\bigoplus_{\mathds{N}}K_{0}(C^*(G)) & \xrightarrow{1-K_0(\tau)} & \bigoplus_{\mathds{N}}K_{0}(C^*(G)) & \rightarrow & K_0(\mathds{U}[\varphi]) \\
\uparrow                            &                      &                                     &                 & \downarrow \\
K_1(\mathds{U}[\varphi])          & \leftarrow           & \bigoplus_{\mathds{N}}K_{1}(C^*(G)) & \xleftarrow{1-K_1(\tau)}& \bigoplus_{\mathds{N}} K_{1}(C^*(G))\\
\end{array}
$$

where $\tau_n(u_g)=s^nu_g{s^*}^n$. Since $K_{0}(\mathbf{K})$ is described only by matrices of the type $E_{i,i}$, consider some $u_{g_i}s^n{s^*}^nu_{g_i^{-1}}u_{\varphi^n(t)}\in D[\varphi]\rtimes_\omega G$. Then\footnote{The projections and unitaries of $D[\varphi]\rtimes_\omega G$ are combinations of elements of that type.}
$$
K_*(\tau)[u_{g_i}s^n{s^*}^nu_{g_i^{-1}}u_{\varphi^n(t)}]_*=[u_{\varphi(g_i)}s^{n+1}{s^*}^{n+1}u_{\varphi(g_i^{-1})}u_{\varphi^{n+1}(t)}]_*,
$$

which implies that $K_{*}(\tau)$ corresponds to a shift in $\bigoplus_{\mathds{N}} K_{*}(C^*(G))$. So denote by $\sigma$ the shift operator, to see that the six-term sequence above turns into
$$
\begin{array}{ccccc}
\bigoplus_{\mathds{N}}K_{0}(C^*(G)) & \xrightarrow{1-\sigma} & \bigoplus_{\mathds{N}}K_{0}(C^*(G)) & \rightarrow & K_0(\mathds{U}[\varphi]) \\
\uparrow                            &                        &                                     &             & \downarrow \\
K_1(\mathds{U}[\varphi])          & \leftarrow             & \bigoplus_{\mathds{N}}K_{1}(C^*(G)) & \xleftarrow{1- \sigma}& \bigoplus_{\mathds{N}} K_{1}(C^*(G))\\
\end{array}
$$

But the application $1- \sigma$ has null kernel and Im$(1-\sigma)$ only contains vectors\nl $(x_0,x_1,\ldots x_n,0,0,\ldots)$ whose sum of coordinates equals zero. This together with the direct limit description (\ref{eq11}) implies that
$$
\dfrac{\bigoplus_{\mathds{N}}K_{*}(C^*(G))}{\hbox{Im }(1-\sigma)}\cong K_{*}(C^*(G))
$$

via
$$
\overline{(x_0,x_1,\ldots x_n,0,0,\ldots)}\mapsto\D\sum_{i=0}^n x_i.
$$

Solving the six-term sequence we get
\begin{equation*}
\begin{split}
K_{*}(\mathds{U}[\varphi])\cong K_{*}(C^*(G)).
\end{split}
\end{equation*}

Therefore, we have the following.
\begin{teo}\label{teo2} Consider $\varphi$ an injective endomorphism with infinite cokernel of some discrete countable group $G$, and construct the C$^*$-algebra $\mathds{U}[\varphi]$ as in Proposition \ref{prop5}. Then $K_{*}(\mathds{U}[\varphi])\cong K_{*}(C^*(G))$.
\end{teo}
\begin{flushright}

  $\square$

  \end{flushright}

\begin{teo}\label{teo3} Consider $\varphi$ a pure injective endomorphism with infinite cokernel of some discrete countable amenable group $G$, and construct the C$^*$-algebra $\mathds{U}[\varphi]$ as in Proposition \ref{prop5}. Then it is classifiable by Kirchberg's classification theorem.
\end{teo}
\begin{flushright}

  $\square$

  \end{flushright}

But note that, if we take two diferent pure injective endomorphisms of some discrete countable amenable group $G$, both C$^*$-algebras will be classifiable by Kirchberg's theorem, and in both objects $K_0(\mathds{U}[\varphi])\ni[1]_0\mapsto [1]_0\in K_0(C^*(G))$. Thus they are isomorphic.

\begin{corol} Satisfied the conditions above, for a fixed group $G$ any choice of endomorphism $\varphi$ generates the same C$^*$-algebra $\mathds{U}[\varphi]$.
\end{corol}
\begin{flushright}

  $\square$

  \end{flushright}
	
\subsection{Semigroup C$^*$-algebra description of $\mathds{U}[\varphi]$}

In \cite{Li0} and \cite{Li2} Li introduced and developed the concept of a C$^*$-algebra associated with a semigroup. His definitions are similar to our C$^*$-algebra associated with an endomorphism, and we will prove that when the semigroup is of the form $S=G\rtimes_\varphi\mathds{N}$ i.e, a semidirect product of a group $G$ with $\mathds{N}$ implemented by an injective endomorphism, the C$^*$-algebra of this semigroup can be viewed as the C$^*$-algebra associated with the endomorphism $\varphi$ defined in this section. This isomorphism together with extra restrictions on our initial data will allow us to conclude similar results concerning the K-theory of $\mathds{U}[\varphi]$ as the one obtained in Theorem \ref{teo2}.

The first clue to suggest this isomorphism is that both constructions use a set of isometries indexed by the semigroup to generate the C$^*$-algebras. And the main step in getting the desired isomorphism is to compare the set of projections used in both definitions and, for this purpose, we shall study the sets which index these projections, namely $B'$ in our case (Lemma \ref{lema1}) and the set $\mathcal{J}$ of constructible right ideals in Li's case (before Definition 2.2 \cite{Li2}). Note that both are defined as a certain set of subsets of the given structure, and they are closed with respect to some set operations.

The problem is that here $B'$ is a set of subsets of a group and Li defines $\mathcal{J}$ containing subsets of a semigroup. However the following holds:

\begin{prop} $\mathcal{J}=\{(g,n)S:\; (g,n)\in S\}$.
\end{prop}
\dem One just have to use the fact that sets of the type
$$
(g,n)S\cap(h,m)S
$$

and
$$
(g,n)^{-1}(h,m)S
$$

are both of the form $(k,l)S$ or $\emptyset$.

This result is also proved in \cite{CuEcLi1} Lemma 6.3.3.
\begin{flushright}

  $\square$

  \end{flushright}

The result above will allow us to establish the isomorphism between the algebra $\mathds{U}[\varphi]$ defined in this chapter and the full semigroup C$^*$-algebra C$^*(S)$ defined by Li in Definition 2.2 of \cite{Li2}.

Consider an endomorphism $\varphi$ of a group $G$ with $B$ containing only subgroups of the form $\varphi^k(G)$. By Proposition \ref{prop5} the C$^*$-algebra $\mathds{U}[\varphi]$ is the universal one generated by
\begin{equation*}
\begin{split}
\hbox{unitaries }&\{u_g:g\in G\}\hbox{ and}\\
\hbox{one isometry }&s
\end{split}
\end{equation*}

satisfying
\begin{enumerate}
	\item[(i)] $u_gs^nu_hs^m=u_{g\varphi^n(h)}s^{n+m}$.
\end{enumerate}

\begin{prop}\label{propli1} We have
\begin{equation*}
\mathds{U}[\varphi]\cong C^*(S),
\end{equation*}
with the latter defined as in \cite{Li2}.
\end{prop}
\dem The C$^*$-algebra $C^*(S)$ is generated by isometries $\{v_{(g,n)}:\;(g,n)\in S\}$ and projections $\{e_X:\; X\in\mathcal{J}\}$ with $\mathcal{J}=\{(g,n)S:\;(g,n)\in S\}$ (by the proposition above).

To prove that the isomorphism holds, first note that the unitaries $v_{(g,0)}$ and the isometries $v_{(e,n)}$ satisfy the relation generating $\mathds{U}[\varphi]$ ((i) above), so there exists a $*$-homomorphism
\begin{equation*}
\begin{split}
\Phi: \mathds{U}[\varphi]&\rightarrow C^*(S)\\
                        u_g&\mapsto v_{(g,0)};\\
                        s^n&\mapsto v_{(e,n)}.
\end{split}
\end{equation*}

For the inverse map, consider the set of isometries $\{u_gs^n:\; (g,n)\in S\}$ and the set of projections
$$
\{u_hs^m{s^*}^mu_{h^{-1}}:\hbox{ associated with }(h,m)S\in\mathcal{J}\}.
$$

Some calculations show that these two sets satisfy the 5 conditions generating $C^*(S)$ (ref. \cite{Li2}). By the universality of this C$^*$-algebra there exists a $*$-homomorphism
\begin{equation*}
\begin{split}
\Psi: C^*_s(S)&\rightarrow \mathds{U}[\varphi]\\
     v_{(g,n)}&\mapsto u_gs^n, and\\
  e_{[(h,m)S]}&\mapsto u_hs^m{s^*}^mu_{h^{-1}}.
\end{split}
\end{equation*}

It is easy to see that $\Phi$ and $\Psi$ are inverses of each other.
\begin{flushright}

  $\square$

  \end{flushright}

\begin{corol}Consider $\varphi$ an injective endomorphism with infinite cokernel of some discrete countable group $G$ and the semidirect product semigroup $S=G\rtimes_{\varphi}\mathds{N}$. Then
$$
K_{*}(C^*(S))\cong K_{*}(C^*(G)),
$$

with $C^*(S)$ as defined in \cite{Li2}.
\end{corol}
\begin{flushright}

  $\square$

  \end{flushright}

There are two more C$^*$-algebras associated with a semigroup $S$. The first one is the concrete representation of $S$ called the reduced semigroup C$^*$-algebra of $S$, denoted by $C^*_r(S)$ and defined in Definition 2.1 in \cite{Li2}. It is easy to check that there exists a surjective $*$-homomorphism
$$
\lambda:C^*(S)\rightarrow C^*_r(S).
$$

For the second one, note that the semigroup $S$ can be viewed as a subsemigroup of the group $\overline{S}$ (defined in the beginning of Section 4.1, Chapter 1), and this allows us to define another C$^*$-algebra associated with $S$, namely $C^*_s(S)$ (Definition 3.2 of \cite{Li2}). It has the same generators as $C^*(S)$ with minor additional relations, so that there is a surjective $*$-homomorphism
$$
\pi_s:C^*(S)\rightarrow C^*_s(S).
$$

But remember that if $\varphi$ is pure and $G$ is amenable the C$^*$-algebra $\mathds{U}[\varphi]$ is simple (and purely infinite) by Theorem \ref{teo1}, and thus so is $C^*(S)$. Therefore we have:

\begin{teo}\label{teo4}For semigroups of the form $S=G\rtimes_\varphi\mathds{N}$ with $G$ an amenable discrete countable group and $\varphi$ a pure injective endomorphism of $G$, the C$^*$-algebras $C^*(S)$, $C^*_s(S)$ and $C^*_r(S)$ defined in \cite{Li2} are isomorphic to $\mathds{U}[\varphi]$. By Theorem \ref{teo2}, we also conclude that
$$
K_{*}(C^*(S))\cong K_{*}(C^*(G)).
$$

Moreover by Theorem \ref{teo3}, they are classifiable by Kirchberg's classification theorem \cite{Kirchcla}.
\end{teo}
\begin{flushright}

  $\square$

  \end{flushright}

The theorem above provides another powerful tool to calculate the K-theory of $\mathds{U}[\varphi]$, which agrees with Theorem \ref{teo2}, just using Theorem 6.3.4 of \cite{CuEcLi1}. For this, note that $G$ being amenable implies that $\overline{S}$ also is (in \cite{CuEcLi1} this group is called the enveloping group of $S$) and therefore it satisfies the Baum-Connes conjecture with coefficients, and the following result applies.

\begin{teo}\label{teo5}For an amenable group $G$ and a pure injective endomorphism with infinite cokernel $\varphi$ of $G$ consider the semigroup $S=G\rtimes_\varphi\mathds{N}$ and choose $B=\{\varphi^k(G)\}$ for some $k\in\mathds{N}$. Then
$$
K_{*}(\mathds{U}[\varphi])\cong K_{*}(C^*(G)).
$$
\end{teo}
\begin{flushright}

  $\square$

  \end{flushright}

\bibliographystyle{plain}
\bibliography{paper}\par\vspace{\baselineskip}

\begin{thebibliography}{10}

\bibitem{BoEx}
G.~Boava and R.~Exel.
\newblock Partial crossed product description of the {C}$^*$-algebras
  associated with integral domains.
\newblock {\em Preprint}, 2010.

\bibitem{Buss1}
A.~Buss.
\newblock A {C}$^*$-álgebra de um {G}rupo.
\newblock Master thesys, Universidade Federal de Santa Catarina, Santa
  Catarina, 2003.

\bibitem{ChoiEffros}
M.-D. Choi and E.~G. Effros.
\newblock Nuclear {C}$^*$-{A}lgebras and the {A}pproximation {P}roperty.
\newblock {\em Amer. J. Math.}, 100:61--79, 1978.

\bibitem{CliPre}
A.~H. Clifford and G.~B. Preston.
\newblock {\em The {A}lgebraic {T}heory of {S}emigroups}.
\newblock Vol I, Mathematical Surveys, No. 7, Amer. Math. Soc., Providence, RI,
  1996.

\bibitem{Cuntz2}
J.~Cuntz.
\newblock Simple {C}$^*$-algebras generated by isometries.
\newblock {\em Comm. Math. Phys.}, 85:173--188, 1977.

\bibitem{CuntzTopMarkovII}
J.~Cuntz.
\newblock A {C}lass of {C}$^*$-algebras and {T}opological {M}arkov {C}hains
  {II}: {R}educible {C}hains and the {E}xt-functor for {C}*-algebras.
\newblock {\em Inventiones mathematicae}, 63:25--40, 1981.

\bibitem{Cuntz1}
J.~Cuntz.
\newblock {\em C$^*$-algebras associated with the $ax+b$-semigroup over
  $\mathds{N}$}.
\newblock Cortiñas, Guillermo (ed.) et al., K-theory and noncommutative
  geometry. Proceedings of the ICM 2006 satellite conference, Valladolid,
  Spain, August 31-September 6, 2006. Zürich: European Mathematical Society
  (EMS). Series of Congress Reports, 2008.

\bibitem{CuEcLi1}
J.~Cuntz, S.~Echterhoff, and X.~Li.
\newblock On the {K}-theory of crossed products by automorphic semigroup
  actions.
\newblock {\em arXiv:1205.5412v1, preprint}, 2012.

\bibitem{Culi1}
J.~Cuntz and X.~Li.
\newblock The {R}egular {C}$^*$-algebra of an {I}ntegral {D}omain.
\newblock {\em Quanta of maths, Clay Math. Proc., Amer. Math. Soc., Providence,
  RI}, 11:149--170, 2010.

\bibitem{CunVer}
J.~Cuntz and A.~Vershik.
\newblock C$^*$-algebras associated with endomorphisms and polymorphisms of
  compact abelian groups.
\newblock {\em arXiv:1202.5960}, 2012.

\bibitem{ExVi}
R.~Exel and F.~Vieira.
\newblock Actions of inverse semigroups arising from partial actions of groups.
\newblock {\em J. Math. Anal. Appl.}, 363:86--96, 2010.

\bibitem{Hirsh}
I.~Hirshberg.
\newblock On {C}$^*$-algebras associated to certain endomorphisms of discrete
  groups.
\newblock {\em New York J. Math.}, 58:99--109, 2002.

\bibitem{Khoska}
M.~Khoshkam and G.~Skandalis.
\newblock Toeplitz algebras associated with endomorphisms and
  {P}imsner-{V}oiculescu exact sequences.
\newblock {\em Pacific Journal of Mathematics}, 181(2):315--331, 1997.

\bibitem{Kirchcla}
E.~Kirchberg.
\newblock The classification of purely infinite {C}$^*$-algebras using
  kasparov's theory.
\newblock {\em Preprint}, third draft, 1994.

\bibitem{Laca1}
M.~Laca.
\newblock From endomorphisms to automorphisms and back: dilations and full
  corners.
\newblock {\em Journal of the London Mathematical Society}, 61:893--904, 2000.

\bibitem{Li1}
X.~Li.
\newblock Ring {C}$^*$-algebras.
\newblock {\em arXiv:0905.4861, preprint}, 2009.

\bibitem{Li0}
X.~Li.
\newblock Nuclearity of semigroup {C}$^*$-algebras and the connection to
  amenability.
\newblock {\em arXiv:1203.0021v2, preprint}, 2012.

\bibitem{Li2}
X.~Li.
\newblock Semigroup {C}$^*$-algebras and amenability of semigroups.
\newblock {\em arXiv:1105.5539v2, preprint}, 2012.

\bibitem{Pivo1}
M.~Pimsner and D.~Voiculescu.
\newblock Exact sequences of {K}-groups and {E}xt-groups of certain crossed
  product {C}$^*$-algebras.
\newblock {\em J. Operator Theory}, 4:549--574, 1980.

\bibitem{Ror}
M.~R$\phi$rdam.
\newblock {\em Classification of {N}uclear {C}$^*$-{A}lgebras}.
\newblock In Classification of {N}uclear {C}$^*$-{A}lgebras. {E}ntropy in
  {O}perator {A}lgebras, {E}ncyclopaedia of {M}athematical {S}ciences, {V}ol.
  126, {S}pringer-{V}erlag, {B}erlin {H}eidelberg {N}ew {Y}ork, 2002.

\bibitem{Scho}
C.~Schochet.
\newblock Topological {M}ethods for {C}$^*$-{A}lgebras {II}: {G}eometric
  {R}esolutions and the {K}ünneth {F}ormula.
\newblock {\em Pacific Journal of Mathematics}, 98:443--458, 1982.

\bibitem{tutu}
J.~L. Tu.
\newblock La conjecture de {B}aum-{C}onnes pour les feuilletages moyennable.
\newblock {\em K-Theory}, 17:215--264, 1999.

\bibitem{vNeu}
J.~von Neumann.
\newblock Zur allgemeinen {T}heorie des {M}asses.
\newblock {\em Fund. {M}ath.}, 13:73--116, 1929.

\end{thebibliography}
\footnotesize DEPARTAMENTO DE MATEM\'{A}TICA, UNIVERSIDADE FEDERAL DE SANTA CATARINA; 88040-900 FLORIAN\'{O}POLIS BRAZIL (felipemate@gmail.com)
\end{document}